\documentclass[12pt]{article}

\usepackage{amsfonts}

\newtheorem{theorem}{Theorem}
\newtheorem{example}{Example}
\newtheorem{lemma}{Lemma}
\newtheorem{corollary}{Corollary}

\newtheorem{proposition}{Proposition}
\newtheorem{assumption}{Assumption}
\newtheorem{algorithm}{Algorithm}
\newtheorem{definition}{Definition}

\def\squareforqed{\hbox{\rlap{$\sqcap$}$\sqcup$}}
\def\qed{\ifmmode\else\unskip\quad\fi\squareforqed}

\def\ba{\begin{array}}
\def\ea{\end{array}}
\def\beann{\begin{eqnarray*}}
\def\eeann{\end{eqnarray*}}
\def\bea{\begin{eqnarray}}
\def\eea{\end{eqnarray}}
\def\beq{\begin{equation}}
\def\eeq{\end{equation}}
\def\BT{\begin{theorem}}
\def\ET{\end{theorem}}
\def\BE{\begin{example}}
\def\EE{\end{example}}
\def\BL{\begin{lemma}}
\def\EL{\end{lemma}}
\def\BP{\begin{proposition}}
\def\EP{\end{proposition}}
\def\BC{\begin{corollary}}
\def\EC{\end{corollary}}
\def\BD{\begin{definition}}
\def\ED{\end{definition}}
\def\BA{\begin{assumption}}
\def\EA{\end{assumption}}
\def\BAL{\begin{algorithm}}
\def\EAL{\end{algorithm}}

\def\R{\mathbb{R}}
\def\cF{{\cal F}}
\def\cG{{\cal G}}
\def\hn{\lfloor n/2 \rfloor}

\begin{document}

\title{An improved Kalai-Kleitman bound for the diameter of a polyhedron}
\author{
Michael J.~Todd
\thanks{
School of Operations Research and Information Engineering,
Cornell University, Ithaca, NY 14853, USA.
E-mail {\tt mjt7@cornell.edu}.
}}
\maketitle
\abstract{
Kalai and Kleitman \cite{KK} established the bound $n^{\log(d) + 2}$ for the 
diameter of a $d$-dimensional polyhedron with $n$ facets. Here we improve the
bound slightly to $(n-d)^{\log(d)}$.
}

\section{Introduction}
A $d$-polyhedron $P$ is a $d$-dimensional set in $\R^d$ that is the 
intersection of a finite number of half-spaces of the form
$H := \{x \in \R^d: a^Tx \leq \beta\}$. 
If $P$ can be written as the intersection of $n$ half-spaces $H_i$, $i = 
1,\dots,n$, but not fewer, we say it has $n$ facets and these facets are the 
faces $F_i = P \cap H_i$, $i = 1,\dots,n$, each affinely isomorphic to a
$(d-1)$-polyhedron with at most $n - 1$ facets. We then call $P$ a 
$(d,n)$-polyhedron.

We say $v \in P$ is a {\em vertex} of $P$ if there is a half-space $H$
with $P \cap H = \{v\}$. (A polyhedron is {\em pointed} if it has a
vertex, or equivalently, if it contains no line.)
Two vertices $v$ and $w$ of $P$ are {\em adjacent}
(and the set $[v,w] := \{(1 - \lambda)v + \lambda w: 0 \leq \lambda \leq 1 \}$
an {\em edge} of $P$) if there is a half-space $H$ with $P \cap H = [v,w]$.
A {\em path of length} $k$ from vertex $v$ to vertex $w$ in $P$ is a
sequence $v = v_0,v_1,\dots,v_k = w$ of vertices with $v_{i-1}$ and $v_i$
adjacent for $i=1,\dots,k$. The {\em distance} from $v$ to $w$ is the
length of the shortest such path and is denoted $\rho_P(v,w)$, and the diameter
of $P$ is the largest such distance,
\[
\delta(P) := \max\{ \rho_P(v,w): v \mbox{ and } w \mbox{ vertices of }P\}.
\]
We define 
\[
\Delta(d,n) := \max\{\delta(P): P \mbox{ a $(d,n)$-polyhedron}\}
\]
and seek an upper bound on $\Delta(d,n)$. It is not hard to see that
$\Delta(d,\cdot)$ is monotonically non-decreasing. Also, the maximum
above can be attained by a {\em simple} polyhedron, one where each vertex
lies in exactly $d$ facets. See, e.g., Klee and Kleinschmidt \cite{KK2}
or Ziegler \cite{Z}. A related paper, Ziegler \cite{Z2}, gives the
history of the Hirsch conjecture that $\Delta_b(d,n) \leq n - d$, where
$\Delta_b(d,n)$ is defined as above but for bounded polyhedra.

In the last few years, there has been an explosion of papers related to the 
diameters of polyhedra and related set systems. Santos \cite{Sa1}
found a counterexample to the Hirsch conjecture, later refined by
Matschke, Santos, and Weibel \cite{MSW}. Eisenbrand, H\"{a}hnle,
Razborov, and Rothkoss \cite{EHRR} showed that a slightly improved
Kalai-Kleitman bound, $n^{\log(d) + 1}$, held for a very general
class of set families abstracting properties of the vertices of $d$-polyhedra
with $n$ facets, which included the ultraconnected set families considered
earlier by Kalai \cite{Ka1}. (This improved bound, for polyhedra, was
presented first in Kalai \cite{Ka2}.)
Another class of set families was introduced
by Kim \cite{Ki}; adding various properties gave families for which this
bound held, or other families where the maximum diameter grew exponentially.
The latter result is due to Santos \cite{Sa2}. Earlier combinatorial
abstractions of polytopes include the abstract polytopes of Adler and
Dantzig \cite{AD} (these satisfy the Hirsch conjecture for 
$n - d \leq 5$) and the duoids of \cite{To1,To2} (these have a lower bound
on their diameter growing quadratically with $n - d$). We also mention the
nice overview articles of Kim and Santos \cite{KS} (pre-counterexample)
and De Loera \cite{DL} (post-counterexample).

Our bound $(n - d)^{\log(d)}$ fits better with the Hirsch conjecture and
is tight for dimensions 1 and 2. Also, more importantly, it is invariant
under linear programming duality. A pointed $d$-polyhedron with $n$
facets can be written as $\{ x \in \R^d: A x \leq b\}$ for some $n \times
d$ matrix $A$ of full rank and some $n$-vector $b$. Choosing an objective
function $c^Tx$ for $c \in \R^d$ gives the linear programming problem
$\max\{c^Tx: A x \leq b\}$, whose dual is $\min\{b^Ty: A^Ty = c, \, y \geq 0\}$.
The feasible region for the latter is affinely isomorphic to a 
pointed polyhedron of dimension at most $n - d$ with at most $n$ facets, and 
equality is possible.  Hence duality switches the dimensions $d$ and $n - d$.

\section{Result}
We prove \hfill

\BT
For $1 \leq d \leq n$, $\Delta(d,n) \leq (n - d)^{\log(d)}$, with 
$\Delta(1,1) = 0$.
\ET

\noindent
(All logarithms are to base 2; note that $(n - d)^{\log(d)} = d^{\log(n - d)}$
as both have logarithm $\log(d) \cdot \log(n - d)$. 
We use this in the proof below.)

The key lemma is due to Kalai and Kleitman \cite{KK}, and was used by them
to prove the bound $n^{\log(d) + 2}$. We give the proof for completeness.
\BL
For $2 \leq d \leq \hn$, where $\hn$ is the largest integer at most $n/2$,
\[
\Delta(d,n) \leq \Delta(d-1,n-1) + 2 \Delta(d,\hn) + 2.
\]
\EL
\begin{proof}
Let $P$ be a simple $(d,n)$-polyhedron and $v$ and $w$ two vertices of $P$
with $\delta_P(v,w) = \Delta(d,n)$. We show there is a path in $P$ from
$v$ to $w$ of length at most the right-hand side above.
If $v$ and $w$ both lie on the same facet, say $F$, of $P$, then
since $F$ is affinely isomorphic to a $(d-1,m)$-polyhedron with
$m \leq n-1$, we have $\rho_P(v,w) \leq \rho_F(v,w) \leq \Delta(d-1,m)
\leq \Delta(d-1,n-1)$ and we are done.

Otherwise, let $k_v$ be the largest $k$ so that there is a set $\cF_v$
of at most $\hn$ facets with all paths of length $k$ from $v$ meeting
only facets in $\cF_v$. This exists since all paths of length 0 meet only
$d$ facets (those containing $v$), whereas paths of length $\delta(P)$
can meet all $n$ facets of $P$. Define $k_w$ and $\cF_{w}$ similarly.
We claim that $k_v \leq \Delta(d,\hn)$ and similarly for $k_w$. Indeed,
let $P_v \supseteq P$ be the $(d,m_v)$-polyhedron ($m_v = | \cF_{v} |
\leq \hn$) defined by just those linear inequalities corresponding to the
facets in $\cF_{v}$. Consider any vertex $t$ of $P$ a distance $k_v$
from $v$, so there is a shortest path from $v$ to $t$
of length $k_v$ meeting only facets in $\cF_{v}$.
But this is also a shortest path in $P_v$, since
if there were a shorter path, it could not be a path in $P$, and thus
must meet a facet not in $\cF_v$, a contradiction. So
\[
k_v = \delta_{P_v}(v,t) \leq \Delta(d,m_v) \leq \Delta(d,\hn).
\]

Now consider the set $\cG_v$ of facets that can be reached in at most $k_v+1$
steps from $v$, and similarly $\cG_w$. Since both these sets contain more than
$\hn$ facets, there must be a facet, say $G$, in both of them. Thus there
are vertices $t$ and $u$ in $G$ and paths of length at most $k_v + 1$ from
$v$ to $t$ and of length at most $k_w + 1$ from $w$ to $u$. Then
\beann
\Delta(d,n) & = & \rho_P(v,w) \\
& \leq & \rho_P(v,t) + \rho_G(t,u) + \rho_P(w,u) \\
& \leq & k_v + 1 + \Delta(d-1,n-1) + k_w + 1 \\
& \leq & \Delta(d-1,n-1) + 2 \Delta(d,\hn) + 2,
\eeann
since, as above, $G$ is affinely isomorphic to a $(d-1,m)$-polyhedron
with $m \leq n-1$.
\end{proof}
\qed

{\bf Proof of the theorem: \rm }
This is by induction on $d + n$.
The result is trivial for $n = d$, since there can be only one
vertex. Next, the right-hand side gives 1 for $d = 1$ ($n = 2$)
and $n - 2$ for $d = 2$, which are the correct values.
For $d = 3$, it gives $(n - 3)^{\log(3)}$, which is greater
than the correct value $n - 3$ established by Klee \cite{Kl1,Kl2,Kl3}.
(We could make the proof more self-contained by establishing the $d=3$
case from the lemma: a general argument deals with $n \geq 13$, but then
there are seven more special cases to check.)
Below we will give a general inductive step for the case $d \geq 4$,
$n - d \geq 8$. Also, the result clearly holds by induction
if $n < 2d$, since then any two vertices lie on a common facet, so
their distance is at most $\Delta(d-1,n-1)$. The remaining cases are 
$d = 4$, $8 \leq n \leq 11$; $d = 5$,
$10 \leq n \leq 12$; $d = 6$, $12 \leq n \leq 13$; and $d = 7$, $n = 14$.
All these cases can be checked easily using the lemma, the
equation $\Delta(d,d) = 0$, and the equations
$\Delta(5,6) = \Delta(4,5) = \Delta(3,4) = \Delta(2,3) = 1$.

Now we deal with the case $d \geq 4$, $n - d \geq 8$.
For this, $\log(n-d) \geq 3$, so we have
\beann
\Delta(d,n) & \leq & \Delta(d-1,n-1) + 2 \cdot \Delta(d,\hn) + 2 \\
 & \leq & (d-1)^{\log(n-d)} + 2 \cdot d^{\log(n/2 - d)} + 2 \\
 & \leq & \left( \frac{d-1}{d}\right)^{\log(n-d)} d^{\log(n-d)} + 
    2 \cdot d^{\log((n - d)/2)} + 2 \\
 & \leq & \left( \frac{d-1}{d}\right)^3 d^{\log(n-d)} + 
  \frac{2}{d} \cdot d^{\log(n-d)} + 2 \\
 & = & \left(1 - \frac{3}{d}  + 
  \frac{3}{d^2} - \frac{1}{d^3} + \frac{2}{d}\right) d^{\log(n-d)} + 2 \\
 & \leq & \left(1 - \frac{1}{d} + 
  \frac{3}{4d} - \frac{1}{d^3}\right) d^{\log(n-d)} + 2 \\
 & \leq & d^{\log(n-d)} - \frac{1}{4d} \cdot d^{\log(n-d)} - 
  \frac{1}{d^3} \cdot d^{\log(n-d)} + 2 \\
 & \leq & d^{\log(n-d)},
\eeann
since each of the subtracted terms is at least one. This completes the proof.
\qed

{\bf Acknowledgement \rm}
Thanks to G\"{u}nter Ziegler and Francisco Santos
for several helpful comments on a previous draft.

\end{document}